\documentclass{amsart}
\usepackage{graphics}

\newtheorem{theorem}{Theorem}[section]
\newtheorem{lemma}[theorem]{Lemma}

\newtheorem{prop}[theorem]{Proposition}

\theoremstyle{remark}

\numberwithin{equation}{section}

\begin{document}
\newcommand{\beqs}{\begin{equation*}}
\newcommand{\eeqs}{\end{equation*}}
\newcommand{\beq}{\begin{equation}}
\newcommand{\eeq}{\end{equation}}
\newcommand{\bin}[2]{\genfrac{(}{)}{0pt}{0}{#1}{#2}}
\newcommand\mymod[1]{\mbox{ mod}\ {#1}}
\newcommand\smymod[1]{\mbox{\scriptsize mod}\ {#1}}
\newcommand\mylabel[1]{\label{#1}}
\newcommand\eqn[1]{(\ref{eq:#1})}

\title[\tiny On representation of an integer by $X^2+Y^2+Z^2$ 
and the modular equations...]
{On representation of an integer  by $X^2+Y^2+Z^2$ 
\\and the modular equations of degree $3$ and $5$.}

\author{Alexander Berkovich}
\address{Department of Mathematics, University of Florida, Gainesville,
Florida 32611-8105}
\email{alexb@ufl.edu}
\thanks{Research was supported in part by NSA grant H98230-09-1-0051}


\date{December 30th, 2010}


\begin{abstract}

I discuss a variety of results involving $s(n)$, the number of representations of 
$n$ as a sum of three squares.
One of my objectives  is to reveal numerous interesting connections between
the properties of 
this function and  certain modular equations of degree 3 and 5.
In particular, I show that 
\beqs
s(25n)=\left(6-\left(-n|5\right)\right)s(n)-5s\left(\frac{n}{25}\right)
\eeqs
follows easily from the well known Ramanujan modular equation of degree $5$.
Moreover, I establish  new relations between $s(n)$ and $h(n)$, $g(n)$,
the number of representations of $n$ by the  ternary quadratic forms 
$$
2x^2+2y^2+2z^2-yz+zx+xy, \quad
x^2+y^2+3z^2+xy,
$$ respectively. \\
Finally, I propose a remarkable new identity for $s(p^2n)- p s(n)$ with $p$ being an odd prime.
This identity makes nontrivial use of the ternary quadratic forms with discriminants $p^2$, $16p^2$.

\end{abstract}
\dedicatory{There are always flowers for those who want to see them} 

\maketitle

\section{Introduction} \label{sec:intro}
\medskip

Let $(a,b,c,d,e,f)(n)$ denote the number of representations of $n$ by the ternary form
$ax^2+by^2+cz^2+dyz+ezx+fxy$. I will assume that $(a,b,c,d,e,f)(n)=0$, whenever $n \not\in Z$.
Let $s(n)$ denote the number of representations of $n$ by ternary form
$x^2+y^2+z^2$. 
In \cite{HS}, Hirschhorn and Sellers proved in a completely elementary manner that
\beq
s(p^2n)=\left(p+1-\left(-n|p\right)\right)s(n)-ps\left(\frac{n}{p^2}\right),
\mylabel{eq:1.1}
\eeq
when $p=3$. Here $(a|p)$ denotes the Legendre symbol.
It should be pointed out that the authors of \cite{HS} proved \eqn{1.1} for all odd prime numbers $p$
by an appeal to the theory of modular forms.

In Section 2, I will show that \eqn{1.1} with $p=5$ follows easily from the well-known identity for 
$\phi(q)^2-\phi(q^5)^2$ with 
\beq
\phi(q)=\sum_{n=-\infty}^{\infty}q^{n^2}.
\mylabel{eq:1.2}
\eeq
Here and throughout, $q$ is a complex number with $|q|<1$.
I will also provide an elementary proof of the following  
\begin{theorem}\label{t1}
If $n\equiv 1,2\mymod 4$, then
\beq
s(25n)-5s(n) =  4(2,2,2,-1,1,1)(n),
\mylabel{eq:1.3}
\eeq
\end{theorem}
\noindent
and 
\begin{theorem}\label{t2}
If $n\equiv 1,2\mymod 4$, then
\beq
s(9n)-3s(n) =  2(1,1,3,0,0,1)(n).
\mylabel{eq:1.4}
\eeq
\end{theorem}
In Section 5, I will show how to remove the parity restrictions in the above theorems
by proving Theorems $5.2$ and $5.3$. Section 6 contains my new Proposition $6.1$, which generalizes Theorems $1.1$, $1.2$, $5.2$ and $5.3$.
A reader with no vested interest in $q$-series may want to proceed directly to Section 6. 
However, a motivated reader may decide to walk  slowly through the initial sections to experience suffering which will later turn into joy. 

Let me point out that two ternary forms $2x^2+2y^2+2z^2-yz+zx+xy$ and $x^2+y^2+3z^2+xy$ both have class number one.
This implies that these forms are both regular \cite{Dick}, \cite{JKS},\cite{Jones}.
For a recent discussion of the relation between the Ramanujan modular equations
and certain ternary quadratic forms the reader is invited to examine \cite{BeJ}.
And it goes without saying that one should not forget the timeless classic \cite{PB}.

I begin by recalling some standard notations, definitions, and useful formulas.
\beq
(a;q)_\infty:=\prod_{j\ge 0}(1-aq^j),
\mylabel{eq:1.5}
\eeq
and
\beq
E(q):=\prod_{j\ge 1}(1-q^j).
\mylabel{eq:1.6}
\eeq
Note that 
\beq
E(-q)=\frac{E(q^2)^3}{E(q^4) E(q)},
\mylabel{eq:N}
\eeq
Ramanujan's general theta-function $f(a,b)$ is defined by
\beq
f(a,b)=\sum_{n=-\infty}^{\infty}a^{\frac{(n-1)n}{2}}b^{\frac{(n+1)n}{2}},\quad |ab|<1.
\mylabel{eq:1.7}
\eeq
In Ramanujan's notation, the celebrated Jacobi triple product identity takes the shape \cite{Bern}, p.35
\beq
f(a,b)=(-a;ab)_\infty(-b;ab)_\infty(ab;ab)_\infty,\quad |ab|<1.
\mylabel{eq:1.8}
\eeq 
Note that $\phi(q)$ can be interpreted as
\beq
\phi(q)=f(q,q)=\frac{E(q^2)^5}{E(q^4)^2 E(q)^2},
\mylabel{eq:1.9}
\eeq
where the product on the right follows easily from \eqn{1.7}.
We shall also require
\beq
\phi(-q)=\frac{E(q)^2}{E(q^2)}.
\mylabel{eq:E}
\eeq 
Next we define
\beq
\psi(q)=f(q,q^3)=\sum_{n=-\infty}^{\infty}q^{2n^2+n}.
\mylabel{eq:1.10}
\eeq
It is not hard to check that 
\beq
\psi(q)=\frac{1}{2}f(1,q)=\sum_{n\ge 0}q^{\frac{(n+1)n}{2}} = \frac{E(q^2)^2}{E(q)},
\mylabel{eq:1.11}
\eeq
\beq
\sum_{n=-\infty}^\infty q^{(4n+1)^2}=\sum_{n=-\infty}^\infty q^{(4n+3)^2}=q\psi(q^8),
\mylabel{eq:1.11b}
\eeq
and that
\beq
f(q,q^9)f(q^3,q^7) =\frac{E(q^{20})E(q^5)E(q^2)^2}{E(q^4)E(q)},
\mylabel{eq:1.12}
\eeq
\beq
f(q,q^4)f(q^2,q^3) =\frac{E(q^{5})^3 E(q^2)}{E(q^{10})E(q)}.
\mylabel{eq:1.13}
\eeq
The function $f(a,b)$ may be dissected in many different ways. We will use the following trivial dissections \cite {Bern}, p.40, p.49
\beq
\phi(q)=\phi(q^4)+2q\psi(q^8), 
\mylabel{eq:1.14}
\eeq
\beq
\phi(q)=\phi(q^9)+2q f(q^3,q^{15}),
\mylabel{eq:1.15}
\eeq
\beq
\phi(q)=\phi(q^{25})+2q f(q^{15},q^{35})+2q^4 f(q^{5},q^{45}).
\mylabel{eq:1.16}
\eeq
We will also require a special case of  Schr\"oter's formula \cite{Bern}, p. 45
\beq
f(a,b)f(c,d) = f(ac,bd)f(ad,bc)+af\left(\frac{b}{c},ac^2d\right)f\left(\frac{b}{d},acd^2\right),
\mylabel{eq:1.17}
\eeq
provided $ab=cd$.
Setting $a=b=c=d=q$ in \eqn{1.17} we obtain 
\beq
\phi(q)^2=\phi(q^2)^2+4q\psi(q^4)^2.
\mylabel{eq:1.19}
\eeq
Iterating, we find that
\beq
\phi(q)^2=\phi(q^4)^2+4q\psi(q^4)^2+ 4q^2\psi(q^8)^2.
\mylabel{eq:1.20}
\eeq
Next, we set $a=q,b=q^9,c=q^3,d=q^7$ in \eqn{1.17} and square the result. This way we have 
\begin{align}
& f(q,q^9)^2f(q^3,q^7)^2  = f(q^4,q^{16})^2f(q^8,q^{12})^2 \nonumber \\
& + 2qf(q^4,q^{16})f(q^8,q^{12})f(q^6,q^{14})f(q^2,q^{18}) 
+ q^2f(q^6,q^{14})^2f(q^2,q^{18})^2.
\mylabel{eq:1.21}
\end{align}
Finally, we multiply both sides in \eqn{1.21} by 
$$
\frac{E(q^4)\phi(q^5)}{E(q^{20})E(q^{10})^2},
$$
and use \eqn{1.9}, \eqn{1.11}, \eqn{1.12} and \eqn{1.13} to arrive at 
\begin{align}
\phi(q)f(q^2,q^8)f(q^4,q^6) & = \psi(q^4)\phi(q^5)\phi(q^{10}) \nonumber \\
& + 2q\psi(q^2)\psi(q^{10})\phi(q^5)+q^2\psi(q^{20})\phi(q^2)\phi(q^5).
\mylabel{eq:1.22}
\end{align}
This result will come in handy in my proof of \eqn{1.3} with $n\equiv 2\mymod 4$.
To deal with the case $n\equiv 1\mymod 4$ in \eqn{1.3} I will require another identity
\beq
\phi(q)\phi(q^{5})+\sum_{m,n}q^{2m^2+2nm+3n^2}= 2\Pi_1(q),
\mylabel{eq:1.23}
\eeq
where 
\beq
\Pi_1(q)=\frac{E(q^{10})E(q^5)E(q^4) E(q^2)}{E(q^{20})E(q)}.
\mylabel{eq:pi1}
\eeq
This formula was discovered and proven in \cite{BeY}. The proof of \eqn{1.23}, given in \cite{BeY}, used only 
a special case of the Ramanujan $_{1}\psi_{1}$ summation formula \cite{Bern2}, p.64.
Multiplying both sides in \eqn{1.23} by $\psi(q^{10})$ and utilizing \eqn{1.11} and \eqn{1.12} we can rewrite \eqn{1.23} as
\beq
\psi(q^{10})\phi(q)\phi(q^5)+\psi(q^{10})\sum_{m,n}q^{2m^2+2nm+3n^2}=2\psi(q^2)f(q,q^9)f(q^3,q^7).
\mylabel{eq:1.24}
\eeq 
 
\bigskip
\section{The ternary implications of the fundamental modular equation of degree $5$} \label{sec:2}
\medskip

In this section we will make an extensive use of a well-known modular equation of degree $5$ 
\beq
\phi(q)^2-\phi(q^{5})^2=4qf(q,q^9)f(q^3,q^7)
\mylabel{eq:2.1}
\eeq
to prove \eqn{1.1} with $p=5$.
We  note that \eqn{2.1} has an attractive companion
\beq
5\phi(q^5)^2-\phi(q)^2=4\Pi_2(q),
\mylabel{eq:2.1.c}
\eeq
where
\beq
\Pi_2(q)=\frac{E(q^{10})^2E(q^4)E(q)}{E(q^{20})E(q^5)}.
\mylabel{eq:2.1.d}
\eeq
Both \eqn{2.1}, \eqn{2.1.c} are discussed in \cite{Bern}.
We remark that the right hand side of \eqn{2.1} was interpreted in terms 
of so-called self-conjugate $5$-cores in \cite{GKS}.
To proceed further I will need  a sifting operator $S_{t,s}$.
It is defined by its action on power series as follows
\beq
S_{t,s}\sum_{n\ge0}c(n)q^n=\sum_{k\ge 0}c(tk+s)q^k.
\mylabel{eq:2.2}
\eeq
Here $t$, $s$ are integers such that  $0\leq s < t$.
Making use of \eqn{1.16}, we find that
\beq
S_{5,0}\phi(q)^2=\phi(q^5)^2+8qf(q,q^9)f(q^3,q^7).
\mylabel{eq:2.3}
\eeq
And so 
\beq
S_{5,0}(\phi(q)^2-\phi(q^5)^2) = -(\phi(q)^2-\phi(q^5)^2) + 8qf(q,q^9)f(q^3,q^7).
\mylabel{eq:2.4}
\eeq
Employing \eqn{2.1} twice, we see that
\beq
S_{5,0}(qf(q,q^9)f(q^3,q^7)) = qf(q,q^9)f(q^3,q^7).
\mylabel{eq:2.5}
\eeq
Analogously, we can check that 
\beq
S_{5,0}\phi(q)^3 =\phi(q^5)^3 + 24q\phi(q^5)f(q,q^9)f(q^3,q^7),
\mylabel{eq:2.6}
\eeq
and that 
\beq
S_{5,1}\phi(q)^3 = 6f(q^3,q^7)(\phi(q^5)^2 + 4qf(q,q^9)f(q^3,q^7))= 6f(q^3,q^7)\phi(q)^2,
\mylabel{eq:2.7}
\eeq
\beq
S_{5,4}\phi(q)^3 = 6f(q,q^9)(\phi(q^5)^2 + 4qf(q,q^9)f(q^3,q^7))= 6f(q,q^9)\phi(q)^2.
\mylabel{eq:2.8}
\eeq
We note, in passing, that 
thanks to \eqn{1.8}, the right hand side in \eqn{2.7} can be rewritten as an infinite product
\begin{align*}
\sum_{n=0}^{\infty}s(5n+1)q^n = 6\prod_{j=1}^{\infty} & (1-q^{2j})^2(1-q^{10j}) \\
& (1+q^{-1+2j})^4(1+q^{-3+10j})(1+q^{-7+10j}).
\end{align*}
Cooper and Hirschhorn studied the generating functions of subsequences of $s(n)$ that 
could be represented by a single, simple infinite product. For example, \eqn{2.7},\eqn{2.8} and \eqn{4.15} are 
the formulas (3.1), (3.2) and (1.1) in \cite{CoH}. \\
With the aid of \eqn{1.16} we can combine \eqn{2.7} and \eqn{2.8} into a single elegant statement
\beq
S_{5,r}(\phi(q)^3 -3 \phi(q)\phi(q^5)^2)=0,
\mylabel{eq:2.9}
\eeq
where $r =1,4$.
Next, we apply $S_{5,0}$ to both sides of \eqn{2.6} to obtain, with a little help from \eqn{2.5}
\beq
S_{25,0}\phi(q)^3 = \phi(q)^3+ 24q\phi(q)f(q,q^9)f(q^3,q^7).
\mylabel{eq:2.10}
\eeq
Subtracting $5\phi(q)^3$ and making use of \eqn{2.1} again, we deduce that 
\begin{align}
S_{25,0}\phi(q)^3-5\phi(q)^3 & = -4\phi(q)^3+ 6\phi(q)(\phi(q)^2-\phi(q^5)^2) \nonumber\\
& = 2(\phi(q)^3 -3 \phi(q)\phi(q^5)^2).
\mylabel{eq:2.11}
\end{align}
Finally, we apply $S_{5,r}$ with $r=1,4$ to  both sides of \eqn{2.11} to find that 
\beq
S_{125,25r}\phi(q)^3-5S_{5,r}\phi(q)^3=0.
\mylabel{eq:2.12}
\eeq
But it is plain that
\beq
\phi(q)^3=\sum_{n=0}^{\infty}s(n)q^n.
\mylabel{eq:2.13}
\eeq
And so the equation \eqn{2.12} can be interpreted as  
\beq
s(25n)-5s(n)=0,
\mylabel{eq:2.14}
\eeq
when $n\equiv 1,4 \mymod 5$.
Thus, the proof of \eqn{1.1} with $p=5$ and  $n\equiv 1,4 \mymod 5$ is complete.\\
We now turn our attention to the $n\equiv 2,3 \mymod 5$ case. 
Subtracting $2\phi(q)^3$ from the extremes of \eqn{2.11}, 
we end up with the formula
\beq
S_{25,0}\phi(q)^3 -7\phi(q)^3=-6\phi(q)\phi(q^5)^2.
\mylabel{eq:2.15}
\eeq
It is now clear that for $r=2,3$
\beq
S_{5,r}(S_{25,0}\phi(q)^3-7\phi(q)^3)=  -6\phi(q)^2 S_{5,r}\phi(q) = 0,
\mylabel{eq:2.16}
\eeq
where in the last step we took advantage of the dissection formula \eqn{1.16}.
Obviously, \eqn{2.16} is equivalent to 
\beq
s(25n)-7s(n)=0,
\mylabel{eq:2.17}
\eeq
when $n\equiv 2,3 \mymod 5$.
And so we completed the proof of \eqn{1.1} with $p=5$ and  $n\equiv 2,3 \mymod 5$.
All that remains to do is  to take care of the $n\equiv 0 \mymod 5$ case. 
Adding $\phi(q)^3$ to both sides of \eqn{2.15} and applying  $S_{5,0}$ to the result,
we get
\beq
S_{5,0}(S_{25,0}\phi(q)^3-6\phi(q)^3) = S_{5,0}(\phi(q)^3-6\phi(q)\phi(q^5)^2).
\mylabel{eq:2.18}
\eeq
Next, we utilize \eqn{1.16}, \eqn{2.1} and \eqn{2.6} to process the right hand side of \eqn{2.18} as follows
\begin{align*}
 S_{5,0}(\phi(q)^3-6\phi(q)\phi(q^5)^2) & = \phi(q^5)^3 +6\phi(q^5)(\phi(q)^2-\phi(q^5)^2)-6\phi(q^5)\phi(q)^2 \\
 & = -5\phi(q^5)^3.
\end{align*}
Hence, we have shown that 
\beq
S_{125,0}\phi(q)^3 -6S_{5,0}\phi(q)^3  = -5\phi(q^5)^3.
\mylabel{eq:2.19}
\eeq
Consequently,
\beq
s(25n)-6 s(n)=-5s\left(\frac{n}{25}\right),
\mylabel{eq:2.20}
\eeq
when $5|n$.
This concludes our proof of \eqn{1.1} with $p=5$.
 
\bigskip
\section{Proof of Theorem $1.1$ } \label{sec:3}
\medskip

I begin by observing that Theorem $1.1$ is equivalent to the following statement
\beq
S_{100,25r}\phi(q)^3 -5S_{4,r}\phi(q)^3 = 4 S_{4,r}T(q),
\mylabel{eq:3.1}
\eeq
where
\beq
T(q):= \sum_{x,y,z}q^{2x^2+2y^2+2z^2-yz+zx+xy}
\mylabel{eq:L}
\eeq
and $r=1,2$.
It is not hard to verify that 
\beq
S_{4,1}T(q)=6S_{4,1}X(1,q),
\mylabel{eq:3.2}
\eeq
and that 
\beq
S_{4,2}T(q)= 3S_{4,2}(X(0,q)+ X(2,q)).
\mylabel{eq:3.3}
\eeq
Here 
\beq
X(r,q):= \sum_{\substack{x,\\ y\equiv -z\equiv r\smymod 4}}q^{2x^2+2y^2+2z^2-yz+zx+xy}.
\mylabel{eq:3.4}
\eeq
It takes very little effort to check that 
\beq
2x^2+2y^2+2z^2-zy+ zx + xy = 2\left(x+\frac{y+z}{4}\right)^2 +\frac{5}{8}(y+z)^2+\frac{5}{4}(y-z)^2.
\mylabel{eq:3.new}
\eeq
Hence
\begin{align}
X(r,q) & =\sum_{\substack{x,\\ y\equiv -z\equiv r\smymod 4}}q^{2\left(x+\frac{y+z}{4}\right)^2 +10(\frac{y+z}{4})^2+20(\frac{y-z}{4})^2} \nonumber\\
& =\sum_{\substack{u,\\ w\equiv v+\frac{r}{2}\smymod 2}}q^{2u^2+10v^2+20w^2},
\mylabel{eq:3.5}
\end{align}
for $r=0,2$. 
It is now evident that
\beq
X(0,q)+X(2,q)=\sum_{u,v,w}q^{2u^2+10v^2+20w^2}=\phi(q^2)\phi(q^{10})\phi(q^{20}).
\mylabel{eq:3.6}
\eeq
Using this last result in \eqn{3.3}, we find that
\beq
S_{4,2}T(q)=3\phi(q^5)S_{4,2}(\phi(q^2)\phi(q^{10})).
\mylabel{eq:3.7}
\eeq
Recalling \eqn{1.14}, we obtain at once that 
\beq
4S_{4,2}T(q) =24\phi(q^5)(\psi(q^4)\phi(q^{10})+6q^2\phi(q^2)\psi(q^{20})).
\mylabel{eq:3.8}
\eeq
We now consider $X(r,q)$ with $r=1,3$.
\beqs
X(r,q)=\sum_{\substack{u,\\ v\equiv w\smymod 2}}q^{2u^2+10v^2+5(2w+r)^2}.
\eeqs
Recalling \eqn{1.11b}, we get
\beq
X(1,q)=X(3,q)=\sum_{u,v,\tilde w} q^{2n^2+10v^2+5(4\tilde w+1)^2}=q^5\phi(q^2)\phi(q^{10})\psi(q^{40}).
\mylabel{eq:3.9}
\eeq
Using \eqn{1.14}, \eqn{3.2} and \eqn{3.9}, we deduce that 
\beqs
S_{4,1}T(q) = 6q\psi(q^{10})S_{4,0}(\phi(q^2)\phi(q^{10}))
=6q\psi(q^{10})(\phi(q^{2})\phi(q^{10}) + 4q^3\psi(q^{4})\psi(q^{20})).
\eeqs
Also, it is not hard to check that
\begin{align}
\sum_{m,n}q^{2m^2+2nm+3n^2} & = \sum_{m,n}q^{2(m+n)^2+10n^2}+q^3\sum_{m,n}q^{2(m+n+1)(m+n)+10(n+1)n} \nonumber\\
& = \phi(q^{2})\phi(q^{10}) + 4q^3\psi(q^{4})\psi(q^{20}).
\mylabel{eq:3.10}
\end{align}
This implies that
\beq
4S_{4,1}T(q) =24 q \psi(q^{10})\sum_{m,n}q^{2m^2+2nm+3n^2}.
\mylabel{eq:3.11}
\eeq
Next, we employ \eqn{2.11} to get 
\beq
S_{100,25r}\phi(q)^3 -5S_{4,r}\phi(q)^3 = 2S_{4,r}(\phi(q)^3-3\phi(q)\phi(q^5)^2).
\mylabel{eq:3.12}
\eeq
With the aid of \eqn{1.14}, \eqn{1.20}, \eqn{2.1}, \eqn{2.1.c} we verify that
\beq
S_{4,1}(\phi(q)^3-3\phi(q)\phi(q^5)^2) =24q\psi(q^2)f(q,q^9)f(q^3,q^7)- 12q\phi(q)\phi(q^5)\psi(q^{10}),
\mylabel{eq:3.13}
\eeq
\beq
S_{4,2}(\phi(q)^3-3\phi(q)\phi(q^5)^2) = -24q\psi(q^2)\psi(q^5)^2+12\phi(q)f(q^2,q^8)f(q^4,q^6).
\mylabel{eq:3.14}
\eeq
Utilizing these results in \eqn{3.12} we obtain
\beq
S_{100,25}\phi(q)^3 -5S_{4,1}\phi(q)^3 = 48q\psi(q^2)f(q,q^9)f(q^3,q^7)- 24q\phi(q)\phi(q^5)\psi(q^{10}),
\mylabel{eq:3.15}
\eeq
\beq
S_{100,50}\phi(q)^3 -5S_{4,2}\phi(q)^3 =-48q\psi(q^2)\psi(q^5)^2+24\phi(q)f(q^2,q^8)f(q^4,q^6).
\mylabel{eq:3.16}
\eeq
Recalling \eqn{3.11}, we see that  \eqn{3.1} with $r=1$ is equivalent to 
\beqs
2\psi(q^2)f(q,q^9)f(q^3,q^7)- \phi(q)\phi(q^5)\psi(q^{10}) = \psi(q^{10})\sum_{m,n}q^{2m^2+2nm+3n^2},
\eeqs
which is, essentially, \eqn{1.24}.
Analogously, employing \eqn{3.8}, we find that \eqn{3.1} with $r=2$ is equivalent to 
\beqs
-2q\psi(q^2)\psi(q^5)^2+\phi(q)f(q^2,q^8)f(q^4,q^6)=\phi(q^5)\psi(q^4)\phi(q^{10})+6q^2\phi(q^2)\phi(q^5)\psi(q^{20}),
\eeqs
which is, essentially, \eqn{1.22}.
The proof of Theorem $1.1$ is now complete.\\

In Section 5 we will generalize Theorem $1.1$. To this end we need to define 
\beq
Y(r,q):= \sum_{\substack{x,\\ y+z\equiv r\smymod 4}}q^{2x^2+2y^2+2z^2-yz+zx+xy},
\mylabel{eq:3.17}
\eeq
where $r=0,1,2,3$.
Observe that the condition $y+z\equiv r\mymod 4$  allows us to introduce new summation variables $u,v,w$, defined as
$x= w-v$, $y=2u+v+r$, $z=2u-v$.
Using \eqn{3.new}, it is easy to see that 
\beqs
2x^2+2y^2+2z^2-zy+ zx + xy = 2r^2+w(2w+r)+ 5v(v+r)+ 5u(2u+r).
\eeqs
Hence
\beq
Y(0,q)=\phi(q^2)\phi(q^5)\phi(q^{10}),
\mylabel{eq:3.18}
\eeq
\beq
Y(2,q)= 4 q^3\phi(q^5)\psi(q^4)\psi(q^{20}),
\mylabel{eq:3.19}
\eeq
\beq
Y(1,q)=Y(3,q)=2q^2\psi(q)\psi(q^5)\psi(q^{10}).
\mylabel{eq:3.20}
\eeq
Employing  \eqn{3.10},\eqn{3.18},\eqn{3.19},\eqn{3.20}, we derive 
\beq
T(q)= \sum_{r=0}^{3}Y(r,q)=  \phi(q^5)\sum_{m,n}q^{2m^2+2nm+3n^2} + 4q^2\psi(q)\psi(q^5)\psi(q^{10}).
\mylabel{eq:3.21}
\eeq
It is easy to see that 
\beq
\sum_{\substack{x,\\ y+z\equiv 1\smymod 2}}q^{2x^2+2y^2+2z^2-yz+zx+xy}= Y(1,q)+Y(3,q)=4q^2\psi(q)\psi(q^5)\psi(q^{10}),
\mylabel{eq:3.22}
\eeq
and that
\beq
\sum_{\substack{x,\\ y+z\equiv 1\smymod 2}}q^{2x^2+2y^2+2z^2-yz+zx+xy}= 2Z(q),
\mylabel{eq:3.23}
\eeq
where
\beq
Z(q):= \sum_{\substack{x,\\y\equiv 0\smymod 2, \\ z\equiv 1\smymod 2}}q^{2x^2+2y^2+2z^2-yz+zx+xy}.
\mylabel{eq:3.24}
\eeq
It is worthwhile to point out that $Z(q)$ has six equivalent representations. For example, one has 
\beqs
Z(q):= \sum_{\substack{x\equiv 0\smymod 2,\\y\equiv 1\smymod 2, \\ z}}q^{2x^2+2y^2+2z^2-yz+zx+xy}.
\eeqs
From \eqn{3.22},\eqn{3.23} we deduce that
\beq
Z(q)= 2q^2\psi(q)\psi(q^5)\psi(q^{10}).
\mylabel{eq:3.25}
\eeq
We conclude this Section that by proving that
\beq
\sum_{\substack{x+y\equiv 1\smymod 2,\\y\equiv z\smymod 2}}q^{2x^2+2y^2+2z^2-yz+zx+xy}=Z(q).
\mylabel{eq:3.26}
\eeq
Indeed, the left hand side of \eqn{3.26} can be rewritten as 
\beqs
\sum_{\substack{x\equiv 0\smymod 2,\\y\equiv 1\smymod 2, \\ z\equiv 1\smymod 2}}q^{2x^2+2y^2+2z^2-yz+zx+xy} + \sum_{\substack{x\equiv 1\smymod 2,\\y\equiv 0\smymod 2, \\ z\equiv 0\smymod 2}}q^{2x^2+2y^2+2z^2-yz+zx+xy}.
\eeqs
Now observe that 
\beqs
\sum_{\substack{x\equiv 1\smymod 2,\\y\equiv 0\smymod 2, \\ z\equiv 0\smymod 2}}q^{2x^2+2y^2+2z^2-yz+zx+xy}=\sum_{\substack{x\equiv 0\smymod 2,\\y\equiv 1\smymod 2, \\ z\equiv 0\smymod 2}}q^{2x^2+2y^2+2z^2-yz+zx+xy}.
\eeqs
And so  the left hand side of \eqn{3.26} becomes
\beqs
\sum_{\substack{x\equiv 0\smymod 2,\\y\equiv 1\smymod 2, \\ z}}q^{2x^2+2y^2+2z^2-yz+zx+xy} = \sum_{\substack{x,\\y\equiv 0\smymod 2, \\ z\equiv 1\smymod 2}}q^{2x^2+2y^2+2z^2-yz+zx+xy}= Z(q),
\eeqs
as desired.
 
\bigskip
\section{Cubic modular identities revisited.} \label{sec:4}
\medskip

As in the last section, I begin by observing that Theorem $1.2$ is equivalent to the following statement
\beq
S_{36,9r}\phi(q)^3 -3S_{4,r}\phi(q)^3 = 4 S_{4,r}\phi(q^3)a(q),
\mylabel{eq:4.1}
\eeq
where
\beqs
a(q):= \sum_{x,y}q^{x^2+xy+y^2},
\eeqs
and $r=1,2$.
The function $a(q)$ was extensively studied in the literature \cite{Bern3}, \cite{Borw}, \cite{BBG}, \cite{HGB}.
It appeared in Borwein's cubic analogue of Jacobi's celebrated theta function identity \cite{Borw}.
I will record below some useful formulas
\beq
4a(q^2)\phi(q^3) = \phi(q)^3+ 3\frac{\phi(q^3)^4}{\phi(q)},
\mylabel{eq:4.2}
\eeq
\beq
a(q) = a(q^3)+ 6q\frac{E(q^9)^3}{E(q^3)},
\mylabel{eq:4.3}
\eeq
\beq
a(q) =  \phi(q)\phi(q^3)+ 4q\psi(q^2)\psi(q^6),
\mylabel{eq:4.4}
\eeq
\beq
a(q) =  2\phi(q)\phi(q^3)-\phi(-q)\phi(-q^3),
\mylabel{eq:4.4.b}
\eeq
\beq
2a(q^2)- a(q) = \frac{\phi(-q)^3}{\phi(-q^3)}
\mylabel{eq:4.4.c}
\eeq
\beq
a(q) =  a(q^4)+ 6q\psi(q^2)\psi(q^6).
\mylabel{eq:4.5}
\eeq
Formula \eqn{4.2} appears as equation (6.4) in \cite{Bern3}.
Identities \eqn{4.3}, \eqn{4.4}, \eqn{4.4.b} and \eqn{4.4.c} are discussed in \cite{BBG}.
In order to prove \eqn{4.5}, the authors of \cite{HGB} have shown that 
\beq
2q\psi(q^2)\psi(q^6) = \sum_{u\not\equiv v\smymod 2}q^{u^2+3v^2}.
\mylabel{eq:4.6}
\eeq
We have at once that
\begin {align}
2q\psi(q^2)\psi(q^6) & = \sum_{\substack{u\equiv 1\smymod 2, \\ v\equiv 0\smymod 2}}
q^{u^2+3v^2}+\sum_{\substack{u\equiv 0\smymod 2, \\ v\equiv 1\smymod 2}}q^{u^2+3v^2} \nonumber \\
& = 2q\psi(q^8)\phi(q^{12})+2q^3\phi(q^4)\psi(q^{24}).
\mylabel{eq:4.7}
\end{align}
Combining  \eqn{4.5} and \eqn{4.7}, we have a pretty neat dissection of $a(q)$ mod $4$
\beq
a(q) =  a(q^4)+ 6q\psi(q^8)\phi(q^{12})+ 6q^3\phi(q^4)\psi(q^{24}).
\mylabel{eq:4.8}
\eeq
In \cite{Shen}, L.C. Shen discussed two well-known modular identities of degree $3$
\beq
\phi(q)^2-\phi(q^3)^2 = 4q\frac{\psi(q)\psi(q^3)\psi(q^6)}{\psi(q^2)},
\mylabel{eq:4.9}
\eeq
and 
\beq
\phi(q)^2+\phi(q^3)^2 = 2\frac{\psi(q) f(q,q^2)f(q^2,q^4)}{\psi(q^2)}.
\mylabel{eq:4.10}
\eeq
Multiplying \eqn{4.9} and \eqn{4.10}, and using 
\beq
f(q,q^2)=\frac{E(q^3)^2E(q^2)}{E(q^6)E(q)},
\mylabel{eq:4.11}
\eeq
\beq
f(q,q^5)=\frac{E(q^{12})E(q^3)E(q^2)^2}{E(q^6)E(q^4)E(q)}
\mylabel{eq:4.12}
\eeq
together with \eqn{1.11} we have 
\beq
\phi(q)^4-\phi(q^3)^4 = 8q\phi(q^3)f(q,q^5)^3.
\mylabel{eq:4.13}
\eeq
Next, we rewrite \eqn{4.13} as 
\beq
\frac{\phi(q)^4}{\phi(q^3)}= \phi(q^3)^3 +8qf(q,q^5)^3.
\mylabel{eq:4.14}
\eeq
Recalling \eqn{1.15}, we can recognize the expression on the right as
\beqs
\phi(q^3)^3 +8qf(q,q^5)^3 =S_{3,0}(\phi(q^9) +2qf(q^3,q^{15}))^3= S_{3,0}\phi(q)^3.
\eeqs
And so 
\beq
S_{3,0}\phi(q)^3 =\frac{\phi(q)^4}{\phi(q^3)}.
\mylabel{eq:4.15}
\eeq
Next, we want to show that 
\beq
S_{9,0}\phi(q)^3 =\frac{4\phi(q)^4-3\phi(q^3)^4}{\phi(q)}.
\mylabel{eq:4.20}
\eeq
To this end, we apply $S_{3,0}$ to  both sides of \eqn{4.15}. Utilizing \eqn{1.15}, we find that
\beq
S_{9,0}\phi(q)^3 =\frac{\phi(q^3)^4 + 4 (8q \phi(q^3)f(q,q^{5})^3)}{\phi(q)}.
\mylabel{eq:4.21}
\eeq
The statement in \eqn{4.20} follows immediately from \eqn{4.13} and \eqn{4.21}.
Moreover, we have 
\beq
S_{9,0}\phi(q)^3-5\phi(q)^3 = -\phi(q)^3 -3\frac{\phi(q^3)^4}{\phi(q)} = -4a(q^2)\phi(q^3),
\mylabel{eq:4.24}
\eeq
where we used \eqn{4.2} in the last step.
Adding $2\phi(q)^3$ to the extremes in \eqn{4.24} we derive 
\beq
S_{9,0}\phi(q)^3-3\phi(q)^3 = 2\phi(q)^3 -4a(q^2)\phi(q^3).
\mylabel{eq:4.25}
\eeq
This result will come in handy in my proof of Theorem $5.2$ in the next section. 

\bigskip
\section{Proof of Theorem $1.2$,Theorem $5.2$ and Theorem $5.3$} \label{sec:5}
\medskip

I begin this section by providing an easy proof of two formulas in \eqn{4.1}.
All I need is the following  
\begin{lemma}\label{5}
If $r=1,2$, then
\beq
S_{4,r}(\phi(q)^3- 2a(q^2)\phi(q^3)) = S_{4,r}(a(q)\phi(q^3)).
\mylabel{eq:5.1}
\eeq
\end{lemma}
Proof: This lemma  is a straightforward corollary of \eqn{1.14}, \eqn{4.5} and \eqn{4.8}.
Next, we apply $S_{4,r}$ with $r=1,2$ to \eqn{4.25} and use \eqn{5.1} to obtain
\beq
S_{36,9r}\phi(q)^3-3S_{4,r}\phi(q)^3 = 2S_{4,r}(\phi(q)^3 -2\phi(q^3)a(q^2))= 2S_{4,r}(a(q)\phi(q^3)),
\mylabel{eq:5.1.a}
\eeq
which is \eqn{4.1}, as desired. The proof of Theorem $1.2$ is now complete.
We can do much better, if we realize that \eqn{5.1} is an immediate consequence of the following elegant result
\beq
\phi(q)^3= \phi(q^3)(a(q)+ 2a(q^2)- 2a(q^4)).
\mylabel{eq:5.2}
\eeq
To prove it, we divide both sides by $\phi(q^3)$ and obtain
\beq
\frac{\phi(q)^3}{\phi(q^3)} =2a(q^2)- a(q)+ 2(a(q)-a(q^4)).
\mylabel{eq:5.3}
\eeq
Using \eqn{4.4.c} and  \eqn{4.5} in \eqn{5.3}, we see that \eqn{5.2} is equivalent to
\beq
\frac{\phi(q)^3}{\phi(q^3)}-\frac{\phi(-q)^3}{\phi(-q^3)} = 12q\psi(q^2)\psi(q^6).
\mylabel{eq:5.4}
\eeq
To verify \eqn{5.4}, I replace $q$ by $-q$ in \eqn{4.4.c} and subtract \eqn{4.4.c} to find with the aid of \eqn{4.4.b}
the following 
\beq
\frac{\phi(q)^3}{\phi(q^3)}-\frac{\phi(-q)^3}{\phi(-q^3)} = a(q)-a(-q)= 3(\phi(q)\phi(q^3)-\phi(-q)\phi(-q^3)).
\mylabel{eq:5.5}
\eeq
Subtracting \eqn{4.4} from \eqn{4.4.b} we obtain 
\beq
\phi(q)\phi(q^3)-\phi(-q)\phi(-q^3) = 4q\psi(q^2)\psi(q^6).
\mylabel{eq:5.6}
\eeq
Hence,
\beq
\frac{\phi(q)^3}{\phi(q^3)}-\frac{\phi(-q)^3}{\phi(-q^3)} = 12q\psi(q^2)\psi(q^6),
\mylabel{eq:5.7}
\eeq
as desired. This completes the proof of \eqn{5.2}.
We are now in a position to improve on \eqn{5.1.a}.
Indeed, it follows from \eqn{4.25} and \eqn{5.2} that
\beq
S_{9,0}\phi(q)^3-3\phi(q)^3 = 2\phi(q^3)a(q)- 4\phi(q^3)a(q^4).
\mylabel{eq:5.8}
\eeq
Consequently, we can extend Theorem $1.2$ as
\begin{theorem}\label{nt2}
\beq
s(9n)-3s(n) = 2(1,1,3,0,0,1)(n)-4(4,3,4,0,4,0)(n).
\mylabel{eq:5.9}
\eeq
\end{theorem}
\noindent
It is worthwhile to point out that Theorem $1.1$ can be extended in a similar manner as 
\begin{theorem}\label{nt1}
\beq
s(25n)-5s(n) =  4(2,2,2,-1,1,1)(n)-8(7,8,8,-4,8,8)(n).
\mylabel{eq:5.10}
\eeq
\end{theorem}
\noindent
It is easy to check  that $(7,8,8,-4,8,8)(n)=0$ when $n\equiv1,2\mymod 4$.
And so \eqn{5.10} reduces to \eqn{1.3}  when $n\equiv 1,2\mymod 4$.
Recalling \eqn{2.11}, we see that all that is required to prove Theorem $5.3$ is
\beq
\phi(q)^3-3\phi(q)\phi(q^5)^2 = 2T(q)- 4\tilde T(q),
\mylabel{eq:5.11}
\eeq
where $T(q)$ was defined in \eqn{L}, and
\beq
\tilde T(q):= \sum_{x,y,z}q^{7x^2+8y^2+8z^2-4yz+8zx+8xy}.
\mylabel{eq:NewL}
\eeq
Making easy changes of summation variables $ y\rightarrow x+y $ and $ z\rightarrow  x+z $ in \eqn{L} we find that
\beq
T(q)= \sum_{x,y,z}q^{7x^2+2y^2+2z^2-yz+4zx+4xy}.
\mylabel{eq:5.12}
\eeq
In a similar fashion one can prove that
\beq
\tilde T(q)= \sum_{\substack{x\equiv y\equiv z\smymod 2}}q^{2x^2+2y^2+2z^2-yz+zx+xy}.
\mylabel{eq:5.13}
\eeq
Combining  \eqn{L}, \eqn{3.23},\eqn{3.25}, \eqn{3.26} and \eqn{5.13}, we can easily derive that
\beq
T(q)-\tilde T(q)= 2Z(q)+ Z(q)=  6q^2\psi(q)\psi(q^5)\psi(q^{10}).
\mylabel{eq:5.14}
\eeq
Hence  we can rewite the right hand side of \eqn{5.11} as 
\beqs
2T(q)- 4\tilde T(q)  = 24q^2\psi(q)\psi(q^5)\psi(q^{10})-2T(q).
\eeqs
Recalling \eqn{3.21}, we  see that  \eqn{5.11} is equivalent to 
\beq
\phi(q)^3-3\phi(q)\phi(q^5)^2 = 16q^2\psi(q)\psi(q^5)\psi(q^{10})-2\phi(q^5)\sum_{m,n}q^{2m^2+2nm+3n^2}.
\mylabel{eq:5.15}
\eeq
To prove the above identity we subtract $2\phi(q)\phi(q^5)^2$ from both sides and use  \eqn{1.23}, \eqn{2.1.c} to find that
\beq
\phi(q)\Pi_2(q)= \phi(q^5)\Pi_1(q)-4q^2\psi(q)\psi(q^5)\psi(q^{10}). 
\mylabel{eq:5.16}
\eeq
Next, we multiply both sides of \eqn{5.16} by
\beqs
\frac{E(q^{20})E(q^5)E(q)}{E(q^{10})^2E(q^4)E(q^2)},
\eeqs 
and use \eqn{E} to end up with 
\beqs
\phi(-q^2)^2-\phi(-q^{10})^2 =  -4q^2\frac{E(q^{20})^3 E(q^2)}{E(q^{10})E(q^4)}.
\eeqs 
Finally, replacing $q^2$  by  $q$ in the above, we  deduce  that \eqn{5.11}  is equivalent to  
\beqs
\phi(-q)^2-\phi(-q^{5})^2 =  -4q \frac{E(q^{10})^3 E(q)}{E(q^{5})E(q^2)}.
\eeqs 
Employing  \eqn{N} and \eqn{1.12}, we see  that the last identity is nothing else but \eqn{2.1} with  $q$  replaced by  $-q$.
Hence \eqn{5.11} is true. This completes my  proof of the Theorem $5.3$.

\bigskip
\section{Bold Proposition } \label{sec:6}
\medskip
 
I now proceed to describe the generalization of Theorem $1.2$ for any odd prime $p$.
Observe that the ternary quadratic form  $x^2+y^2+3z^2+xy$ in this theorem has the discriminants $3^{2}$.
We remind the reader that a discriminant of a ternary form $ax^2+by^2+cz^2+dyz+ezx+fxy$ is defined as 
\beqs
\frac{1}{2}\det
\begin{bmatrix} 2a & f & e \\ f & 2b & d \\ e & d & 2c
\end{bmatrix}.
\eeqs 
Using \cite{Leh} it is easy to check that all ternary forms with the discriminant $p^2$ belong to the same genus, say $TG_{1,p}$.
Let $|\mbox{Aut}(f)|$  denote the number of integral automorphs of a ternary quadratic form $f$, 
and let   $R_{f}(n)$ denote the number of representations of $n$ by $f$.
Let $p$ be an  odd prime and  $n\not\equiv 3\mymod 4$. I propose that 
\beq
s(p^2n)-ps(n) =  48\sum_{f\in TG_{1,p}}\frac{R_{f}(n)}{|\mbox{Aut}(f)|} -96\sum_{f\in TG_{1,p}}\frac{R_{f}\left(\frac{n}{4}\right)}{|\mbox{Aut}(f)|}.
\mylabel{eq:6.16}
\eeq
Clearly, one wants to know if the parity restriction on $n$ in \eqn{6.16} can be removed.
In other words, the question is whether a straightforward generalizion of Theorem $5.2$ exists.
Fortunately, the answer is "yes". However, the answer involves the second genus of ternary forms $TG_{2,p}$
with  discriminant $16p^2$. 
Note that, in general, there are twelve genera of the ternary forms with the discriminant $16p^2$ \cite{Leh}.
However, when  $p\equiv 3\mymod 4$ one can create $TG_{2,p}$ from some binary quadratic form of discriminant $-p$. 
It is a well known fact that all binary forms with the discriminant $-p$ belong to the same genus, say $BG_{p}$.
Let $ax^2+bxz+cz^2$ be some binary form $\in BG_{p}$.
We can convert it into ternary form  
$$
f(x,y,z): = 4ax^2+ py^2+4cz^2+4|b|xz.
$$ 
Next, we extend $f$ to a genus that contains $f$.
This genus is, in fact, $TG_{2,p}$ when  $p\equiv 3\mymod 4$. 
It can be shown that the map 
\beqs
BG_{p}\rightarrow TG_{2,p}
\eeqs
does not depend on which specific binary form from $BG_{p}$ we have choosen as our starting point.
I would like to comment that somewhat similar construction was employed in \cite{BeJ} to define 
the so-called $S$-genus.
Let me illustrate this map for $p = 23$. 
In this case,
\beqs
BG_{23}= \{x^2+xz+6z^2,2x^2+xz+3z^2,2x^2-xz+3z^2\}.
\eeqs
Choosing a binary form $x^2+xz+6z^2$ as a starting point one gets
\begin{align*}
& \{x^2+xz+6z^2\}\rightarrow\{4x^2+23y^2+24z^2+4xz\}\rightarrow \\
& \{4x^2+23y^2+24z^2+4xz,8x^2+23y^2+12z^2+4xz,3x^2+31y^2+31z^2-30yz+2zx+2xy\}.
\end{align*}    
We note that 
\beqs
TG_{2,23}:=\{4x^2+23y^2+24z^2+4xz, 8x^2+23y^2+12z^2+4xz,3x^2+31y^2+31z^2-30yz+2zx+2xy\}
\eeqs  
is just one out of twelve possible genera of the ternary form with the discriminant $8464$.
It is instructive to compare $TG_{2,23}$ and 
\beqs
TG_{1,23}:=\{x^2+6y^2+23z^2+xy, 2x^2+3y^2+23z^2+xy,3x^2+8y^2+8z^2-7yz+2zx+2xy\}.
\eeqs  
Clearly,
\beqs
| TG_{1,23}|= |TG_{2,23}|.
\eeqs
Moreover,
\beqs
|\mbox{Aut}(3x^2+8y^2+8z^2-7yz+2zx+2xy)|= |\mbox{Aut}(3x^2+31y^2+31z^2-30yz+2zx+2xy)|=12,
\eeqs
\beqs
|\mbox{Aut}(x^2+6y^2+23z^2+xy)|= |\mbox{Aut}(4x^2+23y^2+24z^2+4xz)|=8,
\eeqs
\beqs
|\mbox{Aut}(2x^2+3y^2+23z^2+xy)|= |\mbox{Aut}(8x^2+23y^2+12z^2+4xz)|=4.
\eeqs
It is a bit less obvious that 
\beqs
(3,31,31,-30,2,2)(4n)= (3,8,8,-7,2,2)(n),
\eeqs
\beqs
(4,23,24,0,4,0)(4n)= (1,6,23,0,0,1)(n),
\eeqs
\beqs
(8,23,12,0,4,0)(4n)= (2,3,23,0,0,1)(n),
\eeqs
and that 
\beqs
(3,31,31,-30,2,2)(m)=(4,23,24,0,4,0)(m)= (8,12,23,0,0,4)(m)=0,
\eeqs
whenever $m\equiv 1,2\mymod 4$.
I propose that the above properties are, in fact, the signature properties of $TG_{2,p}$.
In other words, for any odd prime $p$ there exists an automorphism preserving bijection 
\beqs
H: TG_{2,p}\rightarrow TG_{1,p},
\eeqs
such that , for any $f\in TG_{2,p}$,
\beqs
|\mbox{Aut}(f)|= |\mbox{Aut} H(f)|,
\eeqs
\beq
R_{f}(4n)= R_{H(f)}(n),
\mylabel{eq:pr}
\eeq
and
\beq
 \hspace{3mm}  R_f(m) = 0 ,\hspace{4mm} \mbox{when}
\hspace{3mm} m \equiv 1,2 \bmod 4. 
\mylabel{eq:pr2}
\eeq
Jagy \cite{WJ} suggested that $TG_{1,p}\cup TG_{2,p}$  does not represent any
integer that is  quadratic residue mod $p$ when $p\equiv 1\mymod 4$,
and when $p\equiv 3\mymod 4$ this union does not represent any
integer that is a quadratic nonresidue mod $p$. That is for any $f\in TG_{1,p}\cup TG_{2,p}$
\beqs
R_{f}(n)= 0,
\eeqs
when $(-n|p)=1$.
In addition, he pointed out that $TG_{2,p}$ represents a proper subset of those numbers represented by $TG_{1,p}$.
Lastly, he observed that both $TG_{1,p}$ and $TG_{2,p}$ are anisotropic at $p$.
I discuss one more example. This time I choose $p=17$.
Here one has
\beqs
TG_{1,17}:=\{3x^2+5y^2+6z^2+yz+2zx+3xy, 3x^2+6y^2+6z^2-5yz+2zx+2xy\},
\eeqs
and
\beqs
TG_{2,17}:=\{7x^2+11y^2+20z^2-8yz+4zx+6xy, 3x^2+23y^2+23z^2-22yz+2zx+2xy\}.
\eeqs
Note that 
\beqs
|\mbox{Aut}(3x^2+5y^2+6z^2+yz+2zx+3xy)|= |\mbox{Aut}(7x^2+11y^2+20z^2-8yz+4zx+6xy)|=4,
\eeqs
\beqs
|\mbox{Aut}(3x^2+6y^2+6z^2-5yz+2zx+2xy)|= |\mbox{Aut}(3x^2+23y^2+23z^2-22yz+2zx+2xy)|=12,
\eeqs
\beqs
(3,23,23,-22,2,2)(4n)=(3,6,6,-5,2,2)(n),
\eeqs
\beqs
(7,11,20,-8,4,6)(4n)=(3,5,6,1,2,3 )(n),
\eeqs
\beqs
(7,11,20,-8,4,6)(m)=(3,23,23,-22,2,2)(m)=0,
\eeqs
whenever $m\equiv 1,2\mymod 4$.
It is worthwhile to point out that there are exactly twelve genera with the discriminant $4624$.
Only three of those have the correct cardinality 
\beqs
|TG_{2,17}|=2,
\eeqs
\beqs
|\{3x^2+6y^2+68z^2+2xy,\quad 10x^2+11y^2+14z^2+2yz+4zx+10xy\}|=2,
\eeqs
\beqs
|\{5x^2+7y^2+34z^2+2xy,\quad 6x^2+12y^2+17z^2+4xy\}|=2.
\eeqs
Note, however, that 
\beqs
|\mbox{Aut}(3x^2+6y^2+68z^2+2xy)|= |\mbox{Aut}(10x^2+11y^2+14z^2+2yz+4zx+10xy)|= 4,
\eeqs
and
\beqs
|\mbox{Aut}(5x^2+7y^2+34z^2+2xy)|= |\mbox{Aut}(6x^2+12y^2+17z^2+4xy)|= 4.
\eeqs
And so, $TG_{2,17}$ is a unique genus with the desired properties. 

I would like to conclude this discussion of $TG_{2,p}$ by providing a more 
explicit description valid in three special cases.
If $p\equiv 3\mymod 4$, then $TG_{2,p}$ is the genus that contains
\beqs
4x^2+py^2+(p+1)z^2+4zx.
\eeqs
I remark that the above form was obtained from the principal binary form $x^2+xz+\frac{p+1}{4}z^2$.
If $p\equiv 2\mymod 3$,  then $TG_{2,p}$ is the genus that contains
\beqs
x^2+\frac{4p+1}{3} y^2 +\frac{4p+1}{3}z^2 + \frac{2-4p}{3}yz+2zx+2xy.
\eeqs
If $p\equiv 5\mymod 8$, then $TG_{2,p}$ is the genus that contains 
\beqs
8x^2+\frac{p+1}{2} y^2 +(p+2)z^2 +2yz+8zx+4xy.
\eeqs
Observe that the smallest prime to escape the above net of three special cases is $p=73$.
I am now ready to unveil the promised extension of \eqn{6.16}.
\begin{prop}\label{p1}
Let $p$ be an odd prime, then
\beq
s(p^2n)-p s(n) =  48\sum_{f\in TG_{1,p}}\frac{R_{f}(n)}{|\mbox{Aut}(f)|} -96\sum_{f\in TG_{2,p}}\frac{R_{f}(n)}{|\mbox{Aut}(f)|}.
\mylabel{eq:6.17}
\eeq
\end{prop} 
The proof of this neat result with $p\ge 7$ is beyond the scope of this paper and will be given in \cite{BeJ2}. 
Note, that \eqn{6.16} follows easily from \eqn{pr},\eqn{pr2} and \eqn{6.17}. \\
Below I illustrate Proposition $6.1$ with some initial examples
\beq
s(7^2n)-7s(n) =  6(1,2,7,0,0,1)(n)-12(4,7,8,0,4,0)(n),
\mylabel{eq:6.41}
\eeq
\begin{align}
s(11^2n)-11s(n) & = 4(3,4,4,-3,2,2)(n)+6(1,3,11,0,0,1)(n) \nonumber \\
& -8(3,15,15-14,2,2)(n)- 12(4,11,12,0,4,0)(n),
\mylabel{eq:6.42}
\end{align}
\beq
s(13^2n)-13s(n) = 12(2,5,5,-3,1,1)(n)-24(8,7,15,2,8,4)(n),
\mylabel{eq:6.43}
\eeq
\begin{align}
s(17^2n)-17s(n) & = 12(3,5,6,1,2,3 )(n)+4(3,6,6,-5,2,2)(n) \nonumber \\
& - 24(7,11,20,-8,4,6)(n)-8(3,23,23,-22,2,2)(n),
\mylabel{eq:6.44}
\end{align}
\begin{align}
s(19^2n)-19s(n) & = 6(1,5,19,0,0,1)(n)+12(4,5,6,5,1,2)(n) \nonumber \\
& - 12(4,19,20,0,4,0)(n)-24(7,11,23,-10,6,2)(n),
\mylabel{eq:6.45}
\end{align}
\begin{align}
s(23^2n)-23s(n) & = 4(3,8,8,-7,2,2)(n)+6(1,6,23,0,0,1)(n) \nonumber \\
& + 12(2,3,23,0,0,1)(n)-8(3,31,31,-30,2,2)(n) \nonumber \\
& - 12(4,23,24,0,4,0)(n)-24(8,23,12,0,4,0)(n),
\mylabel{eq:6.46}
\end{align}
Finally, I note that \eqn{6.41} implies the following impressive identity
\beqs
8q\psi(-q)E(q^2)^2 S_{7,5}(-q;q^2)_\infty = \phi(q)^3+\phi(q^7)\sum_{m,n}(q^{m^2+mn+2n^2}-2q^{4m^2+4mn+8n^2}).
\eeqs

\noindent
\textbf{Acknowledgements}
\medskip

\noindent
I would like to thank Bruce Berndt, Shaun Cooper, Will Jagy, Rainer Schulze-Pillot for their kind interest and helpful discussions.

\bigskip

\bibliographystyle{amsplain}

\end{document}